\documentclass[12pt]{article}
\usepackage{amsfonts}
\usepackage{amscd}

\newcommand{\sect}[1]{\setcounter{equation}{0}\section{#1}}
\newcommand{\subsect}[1]{\subsection{#1}}

\def\be{\begin{equation}}
\def\ee{\end{equation}}
\def\bea{\begin{eqnarray}}
\def\eea{\end{eqnarray}}

\def\1{\'{\i}}
\def\R{{\mathbb R}}
\def\sch{{\cal S}}

\def\jj{{\cal J}_3}
\def\jp{{\cal J}_+}
\def\jm{{\cal J}_-}
\def\ii{{\cal I}}

\def\luisH{{\cal H}}
\def\luisP{{\cal P}}
\def\luisM{{\cal M}}
\def\luisK{{\cal K}}
\def\luisD{{\cal D}}
\def\luisC{{\cal C}}

\begin{document}

\thispagestyle{empty}

 \
\hfill \

\
\vspace{1.25cm}

 \begin{center}

{\Large{\bf{Twist maps for non-standard  quantum algebras }}}

\smallskip

{\Large{\bf{and discrete Schr\"odinger symmetries}}}

\end{center}

\bigskip\bigskip

\begin{center} A. Ballesteros$\dagger$, F.J.  Herranz$\dagger$,
J. Negro$\ddagger$ and L.M. Nieto$\ddagger$
\end{center}

\begin{center} {\it ${\dagger}$ Departamento de F\1sica\\
Universidad de Burgos,   E-09001  Burgos, Spain} \end{center}

\begin{center} {\it ${\ddagger}$ Departamento de F\1sica
Te\'orica\\
Universidad de Valladolid, E-47011  Valladolid, Spain}
\end{center}

\bigskip\bigskip

\begin{abstract}

The minimal twist map 
introduced by Abdesselam {\em et al} \cite{Abde} for the
non-standard (Jordanian) quantum $sl(2,\R)$ algebra is
used to construct the twist maps for two different non-standard
quantum deformations of the (1+1) Schr\"odinger algebra. 
Such deformations are, respectively, 
the symmetry algebras of a space and a
time uniform lattice discretization  of the 
$(1+1)$ free Schr\"odinger equation. It is shown that the
corresponding twist maps connect the usual Lie symmetry
approach to these discrete equations with non-standard quantum
deformations. This
relationship  leads to a clear
interpretation of  the  deformation parameter  as the step of
the uniform (space or time) lattice.

\end{abstract}

\newpage

%%%%%%%%%%%%%%%%%%%%%%%%%%%%%%%%%%%%%%%%%%%%%%

\sect{Introduction}

The non-standard quantum deformation of the $sl(2,\R)$ algebra
(also known as Jordanian or $h$-deformation) 
\cite{Demidov,Zakr,Ohn,nonsb,Ogi,nonsc,nonsd} has induced the
construction of several non-standard quantum
algebras which are naturally related with $sl(2,\R)$ by means
of either a central extension ($gl(2)$) or contraction  (the
$(1+1)$-dimensional Poincar\'e algebra ${\cal P}$).
The contraction of the quantum $gl(2)$ algebra leads to a 
quantum harmonic oscillator $h_4$ algebra which, in turn, can
be interpreted as  a central extension of the Poincar\'e
algebra ${\cal P}$.  These relationships  have been studied in
\cite{boson} within the context of boson representations and 
are displayed in the l.h.s.\ of the following diagram:
$$
\begin{CD}
U_z(sl(2,\R))@>>{\mbox{\footnotesize{central extension}}}>
U_z(gl(2))@>>{\mbox{\footnotesize{Hopf subalgebra}}}>
U_\tau(\sch)\longrightarrow
\mbox{\footnotesize{Discrete time SE}}\\
@VV{\varepsilon \to 0}V
@VV{\varepsilon \to 0}V\\
U_z({\cal P})@>>{\mbox{\footnotesize{central extension}}}>
U_z(h_4)@>>{\mbox{\footnotesize{Hopf subalgebra}}}> 
U_\sigma(\sch)\longrightarrow
\mbox{\footnotesize{Discrete space SE}}
\end{CD}
$$
The cornerstone of  the above four quantum algebras
is the triangular Hopf  algebra 
with generators $J_3,J_+$ verifying $[J_3,J_+]=2 J_+$, and
with classical $r$-matrix  
\cite{Drinfelda,Drinfeldb}
\be
r=z J_3\wedge J_+ .
\label{aaa}
\ee
 The  deformed commutator, coproduct and
universal quantum $R$-matrix   are:
\be 
 [J_3,J_+]=\frac{e^{2 z J_+}-1}{z}\qquad
 \Delta(J_+)=1\otimes J_+ + J_+ \otimes 1\qquad
\Delta(J_3)=1\otimes J_3 + J_3 \otimes e^{2 z J_+}  
\label{aa}
\ee   
\be
{\cal R}=\exp\{-z J_+\otimes J_3\}\exp\{z J_3\otimes J_+\}.
\label{ab}
\ee
This structure is a Hopf (Borel)
subalgebra  of all the  quantum  algebras that we have
previously mentioned. We recall that the  Jordanian twisting
element for the Borel algebra was given in \cite{Ogi},  the
expression for the $R$-matrix (\ref{ab}) was deduced in
\cite{nonsd} and  another construction for ${\cal R}$ can
be found in \cite{nonsc}. The quantum algebra (\ref{aa})
underlies the
approach to physics at
the Planck scale introduced in \cite{Majida, Majidb}.

The triangular nature of a quantum deformation (like
(\ref{aa})) ensures the existence of a twist operator (as the
one given in
\cite{Ogi}) which relates the (classical) cocommutative
coproduct with the (deformed) non-cocommutative one
\cite{Drinfeldb}. This  means that there should exist 
an invertible  twist map which turns the deformed
commutation rules into usual Lie commutators. 
In this respect, a class of twist maps for the non-standard
quantum $sl(2,\R)$ algebra has been explicitly constructed by 
Abdesselam {\em et al} \cite{Abde}; amongst these maps we
will consider the simplest one, which is called the `minimal
twist map'. A first aim of this paper is to implement  the
minimal twist map in the quantum algebras $U_z(gl(2))$,
$U_z({\cal P})$ and
$U_z(h_4)$, showing their relationships with $U_z(sl(2,\R))$
through either contraction or central extension; this task is
carried out in  the section 2. In relation with this kind of
twists, we recall that a different deformation map for
$U_z({\cal P})$ was introduced in \cite{Khorrami} and a
similar construction for
$U_z(e(3))$ was given in
\cite{e3}. Jordanian twists for $sl(n)$ and for some
inhomogeneous Lie algebras have been
also studied in \cite{Kulisha}, and a general construction
of a chain of Jordanian twists has been introduced in 
\cite{Kulishb} and applied to the semisimple Lie algebras of
the Cartan series $A_n$, $B_n$ and $D_n$.

On the other hand,  $gl(2)$ and $h_4$ are Lie
subalgebras of the centrally extended $(1+1)$-dimensional
Schr\"odinger algebra $\sch$ \cite{Hagen,Nied}. These Lie
subalgebra embeddings have been implemented at  a quantum
algebra level and two non-standard quantum
Schr\"odinger algebras have been obtained from them by imposing
that either  $U_z(gl(2))$ or  $U_z(h_4)$ remains as a Hopf
subalgebra. The former \cite{twotime}, $U_\tau(\sch)\supset
U_z(gl(2))$, has been shown to be  the  symmetry algebra of a
time discretization of the heat or (time imaginary)
Schr\"odinger equation (SE) on a uniform lattice. Likewise, 
 as we shall show in this paper,  the latter
\cite{twospace}, $U_\sigma(\sch)\supset U_z(h_4)$, can be
related with a space discretization of the SE also  on a
uniform lattice. These connections are displayed in the
r.h.s.\ of the above diagram. Obviously, when the deformation
parameters $z,\tau,\sigma$ go to zero we recover the usual Lie
algebra picture and the continuous SE. 

In this context of discrete symmetries,   we
recall that  quantum algebras have
been connected with different versions of spacetime lattices
through several algebraic constructions that have no direct
relationship with the usual Lie symmetry theory
\cite{LRN,FS,null}. Recent works \cite{Winternitz} have also
developed new techniques for dealing with the symmetries of
difference or differential-difference equations and have tried
to adapt in this field the standard methods that have been so
successful when applied to differential equations. An
exhaustive study for the discretization on $q$-lattices of
classical linear differential equations has shown that their
symmetries obeyed to $q$-deformed commutation relations with
respect to the Lie algebra structure of the continuous
symmetries
\cite{FVa,FVb,Dobrev}. However, Hopf algebra structures
underlying these $q$-symmetry  algebras have been not found.
When the discretization of linear equations is made on uniform
lattices it seems that the relevant symmetries
preserve the Lie algebra structure \cite{Luismi,Javier}.
Perhaps, this is the reason why  the symmetry approach to
these equations has not been directly related to quantum
algebras.

The second and main objective of  this paper is to relate the
discrete SE's and  their associated
differential-difference symmetry operators provided by  the
quantum Schr\"odinger algebras $U_\sigma(\sch)$, $U_\tau(\sch)$
with the results concerning a discrete SE and operators
obtained in  \cite{Luismi} by following a Lie symmetry
approach. Since the latter operators close the  Schr\"odinger
Lie algebra $\sch$, the crucial point in our procedure is to
find out the twist maps for the Hopf algebras
$U_\sigma(\sch)$, $U_\tau(\sch)$; these   are
straightforwardly obtained from the maps corresponding to
their Hopf subalgebras $U_z(h_4)$ and
$U_z(gl(2))$. These nonlinear changes of basis allow us to
write explicitly the Hopf  Schr\"odinger algebras with
classical commutators and  non-cocommutative coproduct, so
that the connection with the results of \cite{Luismi} can be 
established. The two quantum algebras $U_\sigma(\sch)$,
$U_\tau(\sch)$ are analyzed separately in the sections 3 and
4, respectively. Finally, some remarks close the paper.

%%%%%%%%%%%%%%%%%%%%%%%%%%%%%%%%%%%%%%%%%%%%%%

\sect{Twist maps for non-standard quantum algebras}

%%%%%%%%%%%%%%%%%%%%%%%%%%%%%%%%%%%%%%%%%%%%%%

\subsect{Non-standard quantum $sl(2,\R)$ algebra}

The commutation rules and coproduct of the Hopf algebra 
$U_z(sl(2,\R))$ in the original form deduced in \cite{Ohn} are
given by
\bea
&&[H,X]=2\,\frac{\sinh (zX)}{z}\qquad
[H,Y]=-Y\cosh (z X) - \cosh(z X)\, Y
\qquad
[X,Y]=H\cr
&&\Delta(X)=1\otimes X+X\otimes 1
\qquad
\Delta(Y)=e^{-zX}\otimes Y +Y\otimes e^{zX}\cr
&&\Delta(H)=e^{-zX}\otimes H +H\otimes e^{zX} .
\label{ba}
\eea
The nonlinear invertible map defined by \cite{Negro}
\be
J_+=X\qquad J_3=e^{zX}H\qquad
J_-=e^{zX}\left(Y- z  \sinh(zX) /4 \right)
\label{bb}
\ee
allows us to write the  former structure of $U_z(sl(2,\R))$ as
follows:
\bea 
&&[J_3,J_+ ]= \frac{e^{2 z J_+} -1  } z \qquad 
[J_3,J_-]=-2 J_- +z J_3^2 \qquad  [J_+  ,J_-]= J_3 \cr
\label{eb} 
&&  \Delta (J_+) =1 \otimes J_+  + J_+ \otimes 1 \qquad
\Delta (J_-) = 1 \otimes J_- + J_-\otimes
e^{2 z J_+ } \cr
&& \Delta (J_3) =1 \otimes J_3  + J_3\otimes
e^{2z J_+ } .
\label{bc}
\eea
We remark that in this basis  the universal quantum
$R$-matrix of $U_z(sl(2,\R))$ adopts the factorized expression
(\ref{ab}) and  the classical $r$-matrix is (\ref{aaa}).

Now if we apply to (\ref{bc}) a second  invertible map given by
\be 
\jp=\frac{1-e^{-2 z J_+}}{2z}\qquad
\jj=J_3  \qquad
\jm= J_- -\frac z2 J_3^2  
\label{bd}
\ee
then we find  the classical commutators of $sl(2,\R)$
\be
[\jj,\jp]=2\jp\qquad
[\jj,\jm]=-2\jm\qquad
[\jp,\jm]=\jj
\label{be}
\ee
while the coproduct turns out to be
\bea
&&\Delta(\jp)=1\otimes \jp + \jp\otimes 1 - 2 z \jp \otimes
\jp\cr
&&\Delta(\jj)=1\otimes \jj + \jj\otimes \frac{1}{1- 2z\jp}\cr
&&\Delta(\jm)=1\otimes \jm + \jm\otimes \frac{1}{1- 2z\jp}
- z \jj \otimes \frac{1}{1- 2z\jp}\jj\cr
&&\qquad\qquad
-z^2 (\jj^2+2\jj)\otimes \frac{\jp}{(1- 2z\jp)^2} .
\label{bf}
\eea
Note that although the
new generator $\jp$ is non-primitive, its coproduct satisfies
\be
\Delta((1-2z\jp)^a)=(1-2z\jp)^a\otimes (1-2z\jp)^a
\label{bg}
\ee
since the old generator $J_+$ fulfils  $\Delta(e^{a z
J_+})=e^{a z J_+}\otimes e^{a z J_+}$ for any real number $a$.

The composition of both maps, (\ref{bb}) and (\ref{bd}), gives
rise to 
\be
\jp=\frac{1-T^{-2}}{2z}\qquad
\jj=TH\qquad
\jm=TY-\frac{z}{2} (TH)^2-\frac{z}{8}(T^2-1)
\label{bh}
\ee
 where $T=e^{zX}$. This (invertible) twist map carries the
former Hopf structure of $U_z(sl(2,\R))$ (\ref{ba}) to the
last one characterized by (\ref{be}) and (\ref{bf}) and is the
so called  minimal twist map obtained by Abdesselam {\em et al}
\cite{Abde}.

We  emphasize that  the coproduct (\ref{bf}) has been
explicitly obtained  in a `closed' form which is  worth to
be  compared with the previous literature on nonlinear maps for
the non-standard quantum $sl(2,\R)$
algebra \cite{Abde,Abdea,Abdeb,Aizawa,vander}  since,
in general,   the transformed coproduct has a very complicated
form. In this respect see \cite{Aizawa} where the
corresponding map is used to construct the representation
theory of $U_z(sl(2,\R))$ and also \cite{vander}
where the Clebsch--Gordan coefficients are computed.

%%%%%%%%%%%%%%%%%%%%%%%%%%%%%%%%%%%%%%%%%%%%%%

\subsect{Non-standard quantum $gl(2)$ algebra}

A  non-standard quantum deformation of $gl(2)$ whose
underlying Lie bialgebra is again generated by the classical
$r$-matrix (\ref{aaa}) was constructed in \cite{boson}. The
Hopf algebra $U_z(gl(2))$ reads
\bea 
&&[J_3,J_+ ]= \frac{e^{2 z J_+} -1  } z \qquad 
[J_3,J_-]=-2 J_- +z J_3^2 \cr
&&[J_+  ,J_- ]= J_3 - I\,e^{2z J_+} \qquad [I,\,\cdot\,]=0 
\label{ca} 
\eea 
\be
\begin{array}{l}
\Delta (J_+) =1 \otimes J_+  + J_+ \otimes 1 \qquad
\Delta (J_3) =1 \otimes J_3  + J_3\otimes
e^{2z J_+ }  \cr
\Delta (I) =1 \otimes I  + I \otimes 1 \qquad
 \Delta (J_-) = 1 \otimes J_- + J_-\otimes e^{2 z J_+ } +  z
J_3\otimes I\, e^{2 z J_+ } .
\end{array} 
\label{cb}
\ee 
The  universal quantum $R$-matrix of $U_z(gl(2))$  is also
given by (\ref{ab}).

The twist map which turns (\ref{ca}) into classical commutation
rules is exactly the same as for $U_z(sl(2,\R))$
given in (\ref{bd}) together with $\ii=I$. The resulting
Hopf structure is
\be
[\jj,\jp]=2\jp\qquad
[\jj,\jm]=-2\jm\qquad
[\jp,\jm]=\jj-\ii\qquad [\ii,\,\cdot\,]=0
\label{cc}
\ee
\bea
&&\Delta(\jp)=1\otimes \jp + \jp\otimes 1 - 2 z \jp \otimes
\jp\cr
&&\Delta(\jj)=1\otimes \jj + \jj\otimes \frac{1}{1- 2z\jp}
\cr
&&\Delta(\jm)=1\otimes \jm + \jm\otimes \frac{1}{1- 2z\jp}
- z  \jj \otimes \frac{1}{1- 2z\jp}(\jj-\ii)\cr
&&\qquad\qquad
-z^2 \left( \jj^2+2 \jj\right)\otimes \frac{\jp}{(1-
2z\jp)^2} \cr
&&\Delta(\ii)=1\otimes \ii + \ii\otimes 1.
\label{cd}
\eea
 It is clear that $U_z(gl(2))$ can be seen as an extended 
$U_z(sl(2,\R))$ by $\ii$ which is a central and primitive
generator; if we take $\ii=0$ we recover the results of the
above subsection.

In relation with this construction we recall that
a two-parameter quantum $gl(2)$ algebra, $U_{g,h}(gl(2))$, was
introduced in \cite{Aneva}; it includes 
$U_{z}(gl(2))$ as a particular case  (when both deformation
parameters are identified with $z$).  The Drinfeld twist
operator and map for $U_{g,h}(gl(2))$ were analyzed in
\cite{Aizawab}.

%%%%%%%%%%%%%%%%%%%%%%%%%%%%%%%%%%%%%%%%%%%%%%

\subsect{Non-standard quantum Poincar\'e algebra}

The Hopf algebra   $U_z(sl(2,\R))$ can be contracted to  the
non-standard (1+1) Poincar\'e algebra
by means of the following transformation of the generators and
deformation parameter  \cite{boson}:
\be
P_+=\varepsilon J_+ \qquad P_-=\varepsilon J_- \qquad  K=\frac
12 J_3 \qquad z\to 2\frac {z}{\varepsilon} .
\label{da}
\ee
The limit $\varepsilon\to 0$ gives rise to the Hopf algebra 
  $U_z({\cal P})$ (written in a null-plane basis): 
\be
\begin{array}{l}
\displaystyle{[K,P_+ ]= \frac{e^{ z P_+} -1  } {z} }\qquad 
[K,P_-]=-  P_-  \qquad  [P_+  ,P_- ]=0\cr
  \Delta (P_+) =1 \otimes P_+  + P_+ \otimes 1 
\qquad \Delta (P_-) = 1 \otimes P_- + P_-\otimes
e^{ z P_+ } \cr
\Delta (K) =1 \otimes K  + K\otimes  e^{z P_+ } 
\end{array}  
\label{db}  
\ee
 and the corresponding classical $r$-matrix is $r=  z
 K\wedge P_+$ while the universal $R$-matrix is similar to
(\ref{ab}). 

The contraction  (\ref{da}) allows us to obtain  
straightforwardly   the twist map for $U_z({\cal P})$ and its
resulting Hopf algebra; they are
\be 
 {\cal P}_+=\frac{1-e^{- z P_+}}{z}\qquad
 {\cal K}=K  \qquad
{\cal P}_-= P_-  
\label{dc}
\ee
\be
[{\cal K},{\cal P}_+]= {\cal P}_+\qquad
[{\cal K},{\cal P}_-]=- {\cal P}_-\qquad
[{\cal P}_+,{\cal P}_-]=0 
\label{dd}
\ee
\bea
&&\Delta({\cal P}_+)=1\otimes {\cal P}_+ + {\cal P}_+\otimes 1
-  z {\cal P}_+ \otimes {\cal P}_+\cr
&&\Delta({\cal K})=1\otimes {\cal K} + {\cal K}\otimes
\frac{1}{1-  z{\cal P}_+}\cr &&\Delta({\cal P}_-)=1\otimes
{\cal P}_- + {\cal P}_-\otimes
\frac{1}{1- z{\cal P}_+} .
\label{de}
\eea

%%%%%%%%%%%%%%%%%%%%%%%%%%%%%%%%%%%%%%%%%%%%%%

\subsect{Non-standard quantum harmonic oscillator algebra}

The  Jordanian quantum oscillator algebra
  $U_z( h_4)$ can be obtained by applying the following
contraction  \cite{boson} to the Hopf algebra
$U_z(gl(2))$ given in (\ref{ca})  and (\ref{cb}):
\be
A_+ = \varepsilon J_+ \qquad
A_- = \varepsilon J_- \qquad
N = \frac 12 J_3  \qquad
M = \varepsilon^2 I \qquad
z\to 2\frac {z}{\varepsilon}  
\label{fa}
\ee
together with the limit $\varepsilon \to 0$. Thus we find
 the Hopf structure of $U_z(h_4)$:
\be
\begin{array}{l}
\displaystyle{ [N,A_+]=\frac{e^{z A_+}-1}{z}} \qquad 
[N,A_- ]=-A_-   \qquad [A_- ,A_+]=M  e^{z A_+} \qquad 
[M,\cdot\,]=0 \cr 
\Delta(A_+)=1\otimes A_+ + A_+ \otimes 1 \qquad
 \Delta(N)=1\otimes N +N \otimes e^{z A_+} \cr 
\Delta(M)=1\otimes M +M \otimes 1\qquad \Delta(A_-)=1\otimes
A_-  +A_-  \otimes e^{z A_+}+zN\otimes M\, e^{z A_+}  \cr
\end{array}
\label{fb}
\ee 
whose classical $r$-matrix is $r=z N\wedge A_+$.

The above contraction gives also the  twist map
and the corresponding Hopf algebra:
\be 
 {\cal A}_+=\frac{1-e^{- z A_+}}{z}\qquad
 {\cal N}=N \qquad
 {\cal A}_-= A_-  \qquad {\cal M}=M
\label{fc}
\ee
\be
[{\cal N},{\cal A}_+]= {\cal A}_+\qquad
[{\cal N},{\cal A}_-]=- {\cal A}_-\qquad
[{\cal A}_-,{\cal A}_+]={\cal M}\qquad [{\cal M},\,\cdot\,]=0
\label{fe}
\ee
\bea
&&\Delta({\cal A}_+)=1\otimes {\cal A}_+ + {\cal A}_+\otimes 1
-  z
{\cal A}_+
\otimes
{\cal A}_+\cr
&&\Delta({\cal N})=1\otimes {\cal N} + {\cal N}\otimes
\frac{1}{1- z{\cal A_+}}\cr 
&&\Delta({\cal A}_-)=1\otimes {\cal
A}_- + {\cal A}_-\otimes
\frac{1}{1- z{\cal A}_+} + z {\cal N} \otimes \frac{{\cal
M}}{1- z{\cal A}_+} \cr
 &&\Delta({\cal M})=1\otimes {\cal M} + {\cal M}\otimes 1 .
\label{fg}
\eea

If we set the central generator ${\cal M}=0$ and denote 
${\cal N}={\cal K}$, ${\cal A_+}={\cal P_+}$ and
${\cal A_-}={\cal P_-}$, then we find the results concerning 
 the non-standard (1+1) Poincar\'e algebra of section 2.3.
Therefore we have completed the first part of the diagram of
the Introduction at the level of (minimal) twist maps.

%%%%%%%%%%%%%%%%%%%%%%%%%%%%%%%%%%%%%%%%%%%%%%

\sect{Discrete space Schr\"odinger symmetries}

Let us consider the
discrete version  of  the SE  on a two-dimensional uniform
lattice introduced in
\cite{Luismi}
\be
(\Delta_x^2- 2 m \Delta_t)\phi(x,t)=0  
\label{ac}
\ee
where the  difference operators $\Delta_x$ and $\Delta_t$
  can be expressed in terms of shift
operators 
$T_x= e^{\sigma \partial_x}$ and $T_t=e^{\tau \partial_t}$ as
\be
\Delta_x=\frac{T_x-1}{\sigma} 
\qquad
\Delta_t=\frac{T_t-1}{\tau}. \label{ad} 
\ee
The  parameters $\sigma$ and $\tau$  are the lattice
constants  in the space $x$ and time $t$ directions,
respectively.  The action of $\Delta_x$ (resp.\ $\Delta_t$) on
a function $\phi(x,t)$ consists in a discrete derivative, which
in the limit $\sigma\to 0$ (resp.\  $\tau\to 0$) comes into
$\partial_x$ (resp.\  $\partial_t$). We shall say that an
operator ${\cal O}$ is a symmetry of the linear equation $E
\phi(x,t)=0$ if ${\cal O}$ transforms solutions into 
solutions, that is, if ${\cal O}$ is such that
\be
 E\, {\cal O} = \Lambda\, E  \label{ae}
\ee
 where $\Lambda$ is
another operator. In this way,  the symmetries of the equation
(\ref{ac})  were computed in \cite{Luismi} showing that they
spanned the Schr\"odinger Lie algebra $\sch$, which is
exactly the same result as for the continuous case
\cite{Hagen,Nied}.

%%%%%%%%%%%%%%%%%%%%%%%%%%%%%%%%%%%%%%%%%%%%%%

\subsect{Quantum Schr\"odinger algebra  $U_\sigma(\sch)$ and
discrete space Schr\"odinger equation}

We consider the Schr\"odinger
generators of time translation $H$, space translation $P$,
Galilean boost $K$, dilation $D$, conformal transformation $C$,
and  central generator $M$
\cite{Hagen,Nied}. Let $U_\sigma(\sch)$ be the quantum
Schr\"odinger algebra  \cite{twospace} whose underlying Lie
bialgebra is generated by  the non-standard classical
$r$-matrix
\be
r =\sigma D' \wedge  P  
\label{ga}
\ee
where   $D'=D+\frac 12 M$; hereafter we shall use this notation
in order to simplify  some expressions. The coproduct of
$U_\sigma(\sch)$ has two primitive generators: the central one
$M$ and the space translation
$P$; it reads \cite{twospace}
\be
\begin{array}{l}
\Delta(M)=1\otimes M + M
\otimes 1 \qquad
\Delta(P)=1\otimes P + P \otimes 1  \\  [2pt]
\Delta(H)=1\otimes H + H\otimes e^{2\sigma P} \\  [2pt]
\Delta(K)=1\otimes K + 
K\otimes e^{-\sigma P}+\sigma D'\otimes e^{-\sigma P}M \\ 
[4pt]
\Delta(D)=1\otimes D + D\otimes e^{-\sigma  P} + \frac 12
M\otimes(e^{-\sigma P}-1) \\  [4pt]
\Delta(C)=1\otimes C + C\otimes  e^{-2\sigma P}
-\frac \sigma 2 K\otimes  e^{-\sigma P} D'
+\frac \sigma 2  D'\otimes  e^{-\sigma P}(K-\sigma D' M)  
\end{array} 
\label{gb}
\ee
while the deformed commutation rules are given by
\be
\begin{array}{l}
[D,P]=\frac 1{\sigma} (e^{-\sigma P}-1) \qquad    [D,K]=K
\qquad  [K,P]=M  e^{-\sigma P }\qquad   [M,\,\cdot\,]=0 \\ 
[4pt] [D,H]= - 2 H \qquad\qquad   [D,C]=2 C +\frac {\sigma}2 K
D'  \qquad\qquad [H,P]=0 \\  [4pt]
[H,C]= \frac 12(1+  e^{\sigma P}) D' -\frac 12 M -\sigma K H
\qquad\qquad  [K,C]= \frac {\sigma}2 K^2  \\  [4pt]
[P,C]=-\frac 12 (1+ e^{-\sigma P}) K + \frac {\sigma}2  
e^{-\sigma P} M D'
\qquad\qquad   [K,H]=\frac 1{\sigma}({ e^{\sigma P}-1}) .
\end{array}
\label{gc}
\ee

The relationship between  $U_\sigma(\sch)$ and a space
discretization of the SE can be established by means of 
the following differential-difference realization of (\ref{gc})
in terms of the space and time coordinates $(x,t)$:
\bea
&& P=\partial_x \qquad\qquad H=\partial_t \qquad  \qquad M=
m \cr
&& K=- t \left(\frac {{ e^{\sigma
\partial_x}-1}}{\sigma}\right)   - m x e^{-\sigma
\partial_x}
\qquad D=2  t  \partial_t 
 + x\left( \frac {{  1-e^{-\sigma \partial_x}}}{\sigma}\right)
+\frac 12  \cr
&& C= t^2 \partial_t  e^{\sigma \partial_x} + t x \left(
\frac  {\sinh \sigma\partial_x}{\sigma} 
+\sigma m \partial_t  e^{-\sigma \partial_x}\right)
+\frac 12 m x^2 e^{-\sigma \partial_x}\cr
&& \qquad\qquad -\frac 14 t \left\{  1- 3
e^{\sigma \partial_x} +
  m (1-e^{\sigma \partial_x})\right\}  
-\frac 14 \sigma m (1-   m) x  e^{-\sigma
\partial_x} .
\label{gd}
\eea
The limit $\sigma\to 0$ gives the classical Schr\"odinger 
vector field representation. The Galilei generators
$\{K,H,P,M\}$ close a deformed subalgebra (but not a Hopf
subalgebra) whose Casimir  is
\be
E_{\sigma} = \left(\frac { e^{\sigma P}-1}{\sigma} \right)^2 -
2 M H.
\label{ge}
\ee 
The action of
$E_{\sigma}$ on a function $\phi(x,t)$ through (\ref{gd})
 provides a  space discretization of the SE by choosing for
$E_{\sigma}$ the zero  eigenvalue:
\be
E_{\sigma}\phi(x,t)=0\quad\Longrightarrow\quad \left(    
\left(\frac { e^{\sigma
\partial_x}-1 } {\sigma}\right)^2 - 2 m \partial_t\right)
\phi(x,t)=0 .
\label{gf}
\ee
Furthermore, according to the  definition of a symmetry
operator (\ref{ae}) we find that  the quantum algebra
$U_\sigma(\sch)$ is a symmetry algebra of (\ref{gf}) since
their operators (\ref{gd}) verify
\be
\begin{array}{l}
[E_{\sigma},X]=0\quad \mbox{for}\quad
X\in\{K,H,P,M\}\qquad
[E_{\sigma},D]=  2 E_{\sigma} \cr
 [E_{\sigma},C]=\left\{ t (e^{\sigma \partial_x} +1) +\sigma
m x e^{-\sigma
\partial_x} \right\} E_{\sigma} .
\end{array}
\label{gg}
\ee

%%%%%%%%%%%%%%%%%%%%%%%%%%%%%%%%%%%%%%%%%%%%%%

\subsect{Twist map for $U_\sigma(\sch)$}

The  quantum $h_4$ algebra described in the section 2.4
arises as a Hopf subalgebra of $U_\sigma(\sch)$ once we rename
the generators and deformation parameter of $U_z(h_4)$ as
\be
\begin{array}{l}
A_+=P\qquad A_-=K\qquad N=-D- \frac 12 M \equiv -D'\qquad
z=-\sigma
\end{array}
\label{ha}
\ee
 keeping $M$ as the same central generator. Notice that under
this identification the classical $r$-matrices of both quantum
algebras coincide: $r=zN\wedge A_+=\sigma D'\wedge P$.
The embedding  $U_z(h_4)\subset U_\sigma(\sch)$ allows
us to deduce straightforwardly the (minimal) twist map for
$U_\sigma(\sch)$. The  map associated to $U_z(h_4)$
(\ref{fc}) written in the Schr\"odinger basis (\ref{ha}) reads 
\be 
\luisP=\frac   {e^{\sigma P}-1}{\sigma}
\qquad \luisD=D \qquad \luisK=K \qquad  \luisM=M .
\label{hb}
\ee
The change of basis for $U_\sigma(\sch)$ is completed with the
transformation of the two remaining generators that turns out
to be
\be
\luisH=H\qquad \luisC=C+\frac {\sigma}2 K D' .
\label{hc}
\ee
  In this new basis the commutation
rules (\ref{gc}) of the Hopf algebra $U_\sigma(\sch)$  
come into the Schr\"odinger Lie
algebra:
\be
\begin{array}{llll}
 [\luisD,\luisP]=-\luisP \quad   &[\luisD,\luisK]=\luisK \quad
&[\luisK,\luisP]= \luisM \quad &[\luisM,\,\cdot\,]=0\cr
 [\luisD,\luisH]=-2\luisH \quad  &[\luisD,\luisC]=2\luisC \quad
&[\luisH,\luisC]=\luisD \quad &[\luisH,\luisP]=0 \cr
 [\luisP,\luisC]= - \luisK  \quad  &
 [\luisK,\luisH]= \luisP\quad  &[\luisK,\luisC]=0  
  \quad  &
\end{array}
\label{hd}
\ee
and the coproduct is now given by
\bea
&&\Delta(\luisM)=1\otimes \luisM + \luisM\otimes 1\cr
&&\Delta(\luisP)=1\otimes \luisP + \luisP\otimes 1 + \sigma
\luisP\otimes \luisP\cr
&&\Delta(\luisH)=1\otimes \luisH +
\luisH\otimes (1 + \sigma \luisP)^2\cr
&&\Delta(\luisK)=1\otimes
\luisK + \luisK\otimes  \frac{1}{1 + \sigma \luisP} +\sigma
\luisD'\otimes \frac{\luisM}{1 + \sigma \luisP}\cr
&&\Delta(\luisD)=1\otimes \luisD + \luisD\otimes  \frac{1}{1 +
\sigma \luisP} - 
\frac{1}2 \luisM\otimes \frac{\sigma \luisP}{1 +
\sigma \luisP}\cr
&&\Delta(\luisC)=1\otimes \luisC + \luisC\otimes
\frac{1}{(1 + \sigma \luisP)^2} + \sigma   \luisD'\otimes
\frac{1}{1 + \sigma \luisP}\,\luisK \cr
&&\qquad\qquad+ \frac{\sigma^2}{2}\luisD'(\luisD'-1)\otimes
\frac{\luisM}{(1 +
\sigma \luisP)^2}
\label{he}
\eea
where  we have used again the shorthand notation 
$\luisD'=\luisD+\frac 12 \luisM$. 
Note that  the new 
generator $\luisP$ is  non-primitive but satisfies a property
similar  to  (\ref{bg}).

The mapping defined by (\ref{hb}) and (\ref{hc}) transforms
the  differential-difference realization  (\ref{gd}) into
\bea
&&\luisP=\Delta_x \qquad \qquad \luisH=\partial_t
\qquad \qquad \luisM= m \cr
&&\luisK=- t \Delta_x  -  m   x  
T_x^{-1} \qquad \luisD=2  t  \partial_t  
+ x \Delta_x  T_x^{-1} +
\frac 12  \cr
&&\luisC= t^2 \partial_t   + t x \Delta_x
T_x^{-1}
 +\frac 12 m (x^2 -\sigma x)  T_x^{-2} +\frac 12 t  
\label{hf}
\eea
where $\Delta_x$ and $T_x$ are   the difference 
and shift operators  defined by (\ref{ad}).
Obviously, the Casimir of the Galilei subalgebra $E_\sigma$
(\ref{ge}) leads to the classical one
\be
E=\luisP^2 -2 \luisM\luisH  
\label{hg}
\ee
so that the discretized SE obtained as the realization
(\ref{hf}) of $E\phi(x,t)=0$ is  
\be
(\Delta_x^2- 2 m
\partial_t)\phi(x,t)=0  
\label{hi}
\ee
which coincides with (\ref{gf});  the operators (\ref{hf}) are
symmetries of this equation satisfying
\be
[E,X]=0\quad
X\in\{\luisK,\luisH,\luisP,\luisM\}\qquad [E,\luisD]= 2  E\qquad
 [E,\luisC]= 2 t  E .
\label{hk}
\ee

%%%%%%%%%%%%%%%%%%%%%%%%%%%%%%%%%%%%%%%%%%%%%%

\subsect{Relation of $U_\sigma(\sch)$ with the  Lie symmetry
approach}
 
The discrete space SE  (\ref{hi}) is just the limit $\tau\to
0$ of the  equation (\ref{ac}) considered in \cite{Luismi}.
Hence it is rather natural to expect a connection between the 
differential-difference symmetries obtained in \cite{Luismi}
which close the Lie Schr\"odinger algebra (\ref{hd}) and our
realization of $U_\sigma(\sch)$. Although the operators
(\ref{hf}) do not coincide with those given in  
\cite{Luismi}  we will show that indeed both
realizations are related by means of a similarity
transformation (see \cite{Abde} for 
$U_z(sl(2,\R))$).

The twist map defined by
\bea
&& \luisP=\frac   { e^{\sigma P}-1}{\sigma}
\qquad 
\luisD=D+\frac 12 (1-e^{-\sigma P})
\qquad \luisK= K
-\frac{\sigma}{2} M e^{-\sigma P} \cr
&&\luisM=M \qquad \luisH=H \qquad 
\luisC=C+\frac{\sigma}{2} K D' -\frac{\sigma}{2}
K e^{-\sigma P}-\frac{\sigma^2}8 M e^{-2 \sigma P}   
\label{ia}
\eea
is equivalent to  the one defined by (\ref{hb}) and (\ref{hc})
since it gives rise to the {\it same} Schr\"odinger Lie algebra
(\ref{hd})  and non-cocommutative coproduct (\ref{he}).
However, the new map applied to the  realization (\ref{gd})
leads to
\bea
&&
\luisP=\Delta_x \qquad \qquad \luisH=\partial_t
\qquad \qquad \luisM= m \cr
&& \luisK=- t \Delta_x  -  m   x T_x^{-1}  -\frac{m\sigma}{2}
  T_x^{-1} \qquad \luisD=2  t  \partial_t  
+ x \Delta_x  T_x^{-1}  -
\frac 12 T_x^{-1} +1 \cr
&& \luisC= t^2 \partial_t   + t x \Delta_x
T_x^{-1}
 +\frac 12 m x^2   T_x^{-2}  
 + t \left(  1- \frac 12 T_x^{-1} \right)
-\frac {m \sigma^2}{8}  T_x^{-2}    
\label{ib}
\eea
which are just the symmetry operators of the
equation (\ref{ac}) obtained in \cite{Luismi},  provided the
continuous time limit
$\tau\rightarrow 0$ is performed, $m=\frac 12$  and ${\cal
K}\to -2{\cal K}$.  In other words,   we  have shown  that the
space differential-difference SE introduced  in \cite{Luismi}
has   $U_\sigma(\sch)$ as its  Hopf symmetry
algebra; the operators (\ref{ib}) fulfil the same relations
(\ref{hk}). The deformation parameter
$\sigma$ is interpreted as  the lattice step in the $x$
coordinate, meanwhile the time
$t$ remains a continuous variable.  We also remark that, by
using (\ref{ib}), the solutions of (\ref{hi}) have been
obtained in \cite{Luismi} for $m=\frac 12$.

%%%%%%%%%%%%%%%%%%%%%%%%%%%%%%%%%%%%%%%%%%%%%%

\sect{Discrete time Schr\"odinger symmetries}

%%%%%%%%%%%%%%%%%%%%%%%%%%%%%%%%%%%%%%%%%%%%%%

\subsect{Quantum Schr\"odinger algebra  $U_\tau(\sch)$ and
discrete time Schr\"odinger equation}

A similar procedure can be applied to the quantum Schr\"odinger
algebra $U_\tau(\sch)$ \cite{twotime} coming from the
non-standard classical $r$-matrix
\be
\begin{array}{l}
r= \frac{\tau}{2} D'\wedge H .
\end{array}
\label{ja}
\ee
The coproduct  of $U_\tau(\sch)$ has two primitive generators:
the central one $M$ and the time translation $H$ (instead of
$P$); it is given by \cite{twotime} 
\be
\begin{array}{l}
\Delta(M)=1\otimes M + M
\otimes 1\qquad
\Delta(H)=1\otimes H + H \otimes 1 \\  [2pt]
\Delta(P)=1\otimes P + P\otimes e^{\tau H/2}\\  [2pt]
\Delta(K)=1\otimes K + K\otimes e^{-\tau H/2}+\frac{\tau}{2} 
D'\otimes e^{-\tau H}P\\  [2pt]
\Delta(D)=1\otimes D + D\otimes e^{-\tau H} +  \frac 12
M\otimes (e^{-\tau H}-1)\\  [2pt]
\Delta(C)=1\otimes C + C\otimes  e^{-\tau H}+\frac{\tau}{4}  
D'\otimes e^{-\tau H}  M \end{array}
\label{jb}
\ee
and the compatible commutation rules are
\be
\begin{array}{l}
[D,P]=-P \qquad   [D,K]=K\qquad [K,P]=M   \qquad
[M,\,\cdot\,]=0\\  [2pt]
[D,H]= \frac 2{\tau}
(e^{-\tau H}-1) \qquad   
[D,C]=2 C -\frac{\tau}{2}  (D')^2\qquad [H,P]=0\\ [2pt] 
[H,C]=  D' -\frac 12  M  e^{-\tau H}\qquad
[K,C]=-\frac{\tau}{4}(D'K+KD') \\  [2pt]
[P,C]=-K+\frac{\tau}{4}  (D'P+PD')  \qquad   [K,H]=
e^{-\tau H}P.
\end{array}
\label{jc}
\ee
A differential-difference realization of (\ref{jc})   reads 
\cite{twotime}
\bea
&&H=\partial_t \qquad \qquad
P=\partial_x\qquad\qquad M=m \nonumber\\  [2pt]
&&K=- (t -\tau) e^{-\tau \partial_t}
\partial_x  - m x  \qquad
  D=2 (t-\tau) \left(\frac{1-e^{-\tau
\partial_t} }{\tau}\right) + x\partial_x + \frac 12\nonumber\\ 
[2pt] &&C=(t^2 +\tau b t)
\left(\frac{1-e^{-\tau \partial_t}}{\tau}\right)  + t x
\partial_x  + \frac 12 t  + \frac
12  m x^2   +\tau (b +1) e^{-\tau \partial_t} \nonumber \\ 
[4pt]
&&\qquad\qquad +\frac{\tau}{4} x^2
\partial_{x}^2  +\frac{\tau}{2} (b+1) x \partial_x
+\frac{\tau}{4}(b+1/2)^2  
\label{jd}
\eea
where $b= \frac m2 -2$. A time  discretization of the SE
is obtained by considering the deformed Casimir of the Galilei
subalgebra 
\be
E_{\tau}=P^2- 2 M \left( \frac{e^{\tau H}-1}{\tau}\right)
\label{je}
\ee 
written in terms of the realization (\ref{jd}):
\be
E_{\tau}\phi(x,t)=0\quad\Longrightarrow\quad 
\left(\partial_x^2 - 2 m \left(\frac
{ e^{\tau\partial_t}-1}{\tau}\right)\right)
\phi(x,t)=0 .
\label{jf}
\ee
Under the realization (\ref{jd})  the generators of
$U_\tau(\sch)$ are symmetry operators of this equation
as they satisfy
\be
\begin{array}{l}
[E_{\tau},X]=0\quad X\in\{K,H,P,M\}\qquad [E_{\tau},D]=  2
E_{\tau} \cr
 [E_{\tau},C]=2\left\{ t -\frac{\tau}{4} (1 - m - 2  x
\partial_x)\right\} E_{\tau} .
\end{array}
\label{jg}
\ee

%%%%%%%%%%%%%%%%%%%%%%%%%%%%%%%%%%%%%%%%%%%%%%

\subsect{Twist map for $U_\tau(\sch)$}

The quantum $gl(2)$ algebra studied in the section 2.2 arises
as a Hopf subalgebra of $U_\tau(\sch)$ under the following
identification:
\be
\begin{array}{l}
J_+=H\quad J_-=-C\quad J_3=-D- \frac 12 M \equiv -D'\quad
I=-  \frac 12 M  \quad
z=-\frac 12 \tau .
\end{array}
\label{ka}
\ee
In the Schr\"odinger basis the twist map for
$U_z(gl(2))$ (\ref{bd}) (which is the same as for
$U_z(sl(2,\R))$) is  given by
\be 
\luisH=\frac   {e^{\tau H}-1}{\tau}
\qquad \luisD=D \qquad \luisC=C-\frac{\tau}{4}(D')^2 \qquad 
\luisM=M .
\label{kb}
\ee
The twist map for $U_\tau(\sch)\supset
U_z(gl(2))$ is completed with  the  transformation of the two
generators out of $U_z(gl(2))$ which is simply
\be
\luisP=P\qquad \luisK=K  .
\label{kc}
\ee
If we apply (\ref{kb}) and (\ref{kc}) to (\ref{jb}) and
(\ref{jc}),  we find again the classical commutation rules of
the Schr\"odinger algebra (\ref{hd}) while the coproduct reads
now   
\bea
&&\Delta(\luisM)=1\otimes \luisM +
\luisM\otimes 1\cr
&&\Delta(\luisH)=1\otimes \luisH + \luisH\otimes 1 +
\tau \luisH\otimes \luisH\cr
&&\Delta(\luisP)=1\otimes \luisP +
\luisP\otimes (1 + \tau \luisH)^{1/2}\cr
&&\Delta(\luisK)=1\otimes
\luisK + \luisK\otimes 
\frac{1}{(1 + \tau \luisH)^{1/2}} +\frac{\tau}{2} 
\luisD'\otimes \frac{\luisP}{1 + \tau \luisH}\cr
&&\Delta(\luisD)=1\otimes \luisD
 + \luisD\otimes  \frac{1}{1 + \tau \luisH}
 - \frac{1}2 \luisM\otimes \frac{\tau \luisH}{1 + \tau
\luisH}\cr 
&&\Delta(\luisC)= 1\otimes \luisC + \luisC\otimes 
\frac{1}{1 + \tau
\luisH}   - \frac {\tau}{2} \luisD'\otimes 
  \frac{1}{1 + \tau \luisH}
\luisD \cr
&&\qquad\qquad + \frac{\tau}4 \luisD'(\luisD' - 2 )\otimes
\frac{\tau \luisH}{(1 + \tau \luisH)^2} .
\label{kd}
\eea
In this new basis  the
realization (\ref{jd}) turns out to be
\bea
&& \luisH=\Delta_t \qquad
\luisP=\partial_x  \qquad \luisM= m \cr
&& \luisK=-( t-\tau) \partial_x
T_t^{-1}  - m x \qquad \luisD=2 ( t-\tau)  \Delta_t   T_t^{-1}
+ x \partial_x +  \frac 12  \label{ke}\\
&& \luisC= t^2 \Delta_t   T_t^{-2}    +   x(t-\tau)
\partial_x T_t^{-1}
 +\frac 12 m x^2+3 t T_t^{-2} -\frac 5 2 t T_t^{-1} - 2 \tau
T_t^{-2} +\frac {3\tau}{2} T_t^{-1}  .
\nonumber
\eea
The corresponding
discretized  SE is provided by the (non-deformed)  Casimir $E$
of the Galilei subalgebra (\ref{hg}) written through
(\ref{ke}) leading again to (\ref{jf}): 
\be 
(\partial^2_x- 2 m
\Delta_t)\phi(x,t)=0 .
\label{kf}
\ee
The new operators (\ref{ke}) are symmetries of this equation
satisfying
\be
[E,X]=0\quad
X\in\{\luisK,\luisH,\luisP,\luisM\}\qquad 
[E,\luisD]=  2  E \qquad
[E,\luisC]= 2 (t-\tau) T_t^{-1} E .
\label{kg}
\ee

%%%%%%%%%%%%%%%%%%%%%%%%%%%%%%%%%%%%%%%%%%%%%%

\subsect{Relation of $U_\tau(\sch)$ with the  Lie symmetry
approach}

The connection with the time discretization  of the  SE
analyzed in \cite{Luismi} is provided by the twist map
defined by
\bea
&&
\luisH=\frac   {e^{\tau H}-1}{\tau} \qquad 
\luisD=D+ 2(1-e^{-\tau H}) \qquad 
\luisC=C-\frac{\tau}{4} (D')^2 +\tau D e^{-\tau H} \cr
&&\luisM=M\qquad \luisP= P \qquad \luisK= K -\tau P e^{-\tau
H} .
\label{la}
\eea
This nonlinear map is a similarity transformation of
the former change of basis defined by (\ref{kb}) and (\ref{kc})
since it leads to the same  Lie Schr\"odinger
commutators (\ref{hd}) and non-cocommutative coproduct
(\ref{kd}). Under this  map  the
realization (\ref{jd}) becomes
\bea
&& \luisH=\Delta_t \qquad
\luisP=\partial_x  \qquad \luisM= m \cr
&& \luisK=- t \partial_x
T_t^{-1}  - m x \qquad \luisD=2  t  \Delta_t   T_t^{-1} + x
\partial_x +  \frac 12  \cr
&& \luisC= t^2 \Delta_t   T_t^{-2}    + t x
\partial_x T_t^{-1}
 +\frac 12 m x^2   + t \left( T_t^{-2} 
 - \frac 12 T_t^{-1} \right)  .
\label{lb}
\eea
These difference-differential operators are the limit
$\sigma\rightarrow 0$  of the symmetry operators obtained in
\cite{Luismi} once we set $m=\frac 12$ and
${\cal K}\to -2{\cal K}$. The corresponding discretized  SE is
again  (\ref{kf}) and  the new operators   are
symmetries of this equation satisfying
\be
[E,X]=0\quad
X\in\{\luisK,\luisH,\luisP,\luisM\}\qquad 
[E,\luisD]=  2  E \qquad
[E,\luisC]= 2 t T_t^{-1} E .
\label{lc}
\ee
Henceforth we have explicitly shown that the space
discretization of the SE on a uniform lattice formerly studied
in \cite{Luismi} within a pure Lie algebra approach has
actually a quantum algebra symmetry associated to the Hopf
algebra $U_\tau(\sch)$. Consequently the deformation parameter
$\tau$ is the time  lattice step on this discrete time SE (the
space coordinate $x$ remains as a continuous variable). In this
way, the relationships displayed in the r.h.s. of the diagram
of the Introduction have been studied at the level of twist
maps.

%%%%%%%%%%%%%%%%%%%%%%%%%%%%%%%%%%%%%%%%
\sect{Concluding remarks}

We have explicitly shown that the symmetry algebra
\cite{Luismi} of the space discretization of the SE obtained
from (\ref{ac}) by taking the limit $\tau\to 0$ is just the
quantum Schr\"odinger algebra  $U_\sigma(\sch)$
\cite{twospace} and the deformation parameter $\sigma$ is
exactly  the  space lattice constant.   Likewise, we have also
shown   that the time discretization of the SE
obtained from (\ref{ac}) by means of the limit $\sigma\to 0$
has the quantum Schr\"odinger algebra $U_\tau(\sch)$
\cite{twotime} as its symmetry algebra; in this case, the time
lattice step $\tau$ plays the role of the deformation
parameter. Consequently, a direct relationship between
non-standard (or Jordanian) deformations and regular lattice
discretizations has been established.

We wish to point out that the existence  of a Hopf
algebra structure for the symmetries of a given equation
associated to an elementary system allows us to write
equations of composed systems keeping the same symmetry algebra
\cite{EnricoA,EnricoB}. In order to use this property for the
two discrete SE's here discussed, we see that only the last
commutator in either (\ref{hk}) or (\ref{kg}) involving the
conformal generator  $\luisC$ is not algebraic, but depends
explicitly on the chosen representation (the same happens at
the continuous level). Therefore the composed systems
characterized by the equation $\Delta(E) \phi = 0$ will have,
by construction, $\Delta(\luisH),\Delta(\luisP),
\Delta(\luisK), \Delta(\luisD)$, and $\Delta(\luisM)$ as
symmetry operators  (moreover they close a Hopf subalgebra!).
However, in general,  this will not be the case for
$\Delta(\luisC)$, and a further study on the behaviour of this
operator is needed in order to construct coupled equations
with full quantum Schr\"odinger algebra symmetry.

Finally, we stress that the applicability of the constructive
approach presented here is not limited to the cases analyzed
before, since it could be directly extended to other quantum
algebras by means of their corresponding
differential-difference realizations. In particular, the
results of this paper indicate that there should  exist an
analogous relationship between the discrete symmetries of the
$(1+1)$ wave equation on a uniform lattice obtained in
\cite{Javier} and some non-standard quantum
deformation of the algebra $so(2,2)$. Work on this line is  in
progress.

\newpage

%%%%%%%%%%%%%%%%% ACKNOWLEDGMENTS %%%%%%%%%%%

\noindent {\section*{Acknowledgments}}

\noindent
This work was partially supported by DGES   (Project
PB98--0370) from the Ministerio de Educaci\'on y Cultura  de
Espa\~na and by Junta de Castilla y Le\'on (Projects CO2/197 
 and CO2/399).

\end{document}